\documentclass[]{article}
\usepackage{amssymb}
\usepackage{amsmath}
\usepackage{amssymb}

\newtheorem{thm}{Theorem}

\newtheorem{lem}[thm]{Lemma}

\newtheorem{cj}[thm]{Conjecture}
\newtheorem{col}[thm]{Corollary}

\newcommand{\prf}{\noindent{\it Proof.} }
\newcommand {\cbdo}{\hfill$\Box$}
\def\mod{\mathop{\rm mod}\nolimits}

\begin{document}

\title{A note on $k$-cordial $p$-uniform hypertrees}
\author{Sylwia Cichacz, Agnieszka G\"{o}rlich\\
\small{ University of
Science and Technology AGH, Al. Mickiewicza 30, 30-059 Krak\'ow,
Poland}} \maketitle

\begin{abstract} 
Hovey \cite{Hovey} introduced a $k$-cordial labeling of graphs as a generalization both of harmonious and cordial labelings. He proved that all tress are $k$-cordial for $k \in \{1,...,5\}$ and he conjectured that all trees are
$k$-cordial for all $k$. \\
\indent We consider a corresponding problem for hypergraphs, namely, we show that $p$-uniform hypertrees are $k$-cordial for certain values of $k$.

\end{abstract}

\section{Introduction}
A hypergraph $H$ is a pair $H = (V,E)$ where $V$ is a set of
vertices and $E$ is a set of non-empty subsets of $V$ called
hyperedges. The order of a hypergraph $H$ is denoted by $|H|$ and
the size is denoted by $\|H\|$. If all edges have the same
cardinality $p$, the hypergraph is said to be $p$-uniform. Hence a
graph is $2$-uniform hypergraph. The degree of a vertex $v$, denoted
by $d(v)$, is defined as $d(v) = |{e \in E : v \in e}|$; i.e., the
degree of $v$ is the number of edges to which it belongs. Two
vertices in a
hypergraph are adjacent if there is an edge containing both of them.\\
\indent A {\it walk} in a hypergraph is a sequence $v_0, e_1, v_1,
\ldots, v_{n-1}, e_n, v_n$, where $v_i \in  V$, $e_i \in E$ and
$v_{i-1},v_i \in e_i$ for all $i$. We define a {\it path} in a
hypergraph to be a walk with all $v_i$ distinct and all $e_i$
distinct. A {\it cycle} is a walk containing at least two edges, all
$e_i$ are distinct and all $v_i$ are distinct except $v_0 = v_n$. A
hypergraph is connected if for every pair of its vertices $v,u$,
there is a path starting at $v$ and ending at $u$.
A {\it hypertree} is a connected hypergraph with no cycles.\\
\indent For a $p$-uniform hypergraph $H=(V,E)$ and a $k$-labeling $c: V \to Z_k$ let $v_c(i)=|c^{-1}(i)|$.
The coloring $c$ is said to be {\it $k$-friendly} if $|v_c(i)-v_c(j)| \leq 1$ for any $i \neq j; i,j \in Z_k$.
The coloring $c$ induces an edge labeling $c^*:E \to Z_k$  defined by $c^*(e)=\sum_{v \in e}c(v) \mod k$. Let $e_{c^*}(i)=|{c^*}^{-1}(i)|$.
A hypergraph is said to be {\it $k$-cordial} if it admits such $k$-friendly coloring $c$ that $|e_{c^*}(i)-e_{c^*}(j)|\leq 1$
for any $i \neq j; i,j \in Z_k$. Then we say that the edge coloring $c^*$ is {\it $k$-cordial}. We call a hypergraph \emph{cordial} if it is $2$-cordial.\\
Cordial labeling of graphs was introduced by Cahit \cite{Cahit} as a weakened
version of graceful labelling and harmonious labelling. He proved~\cite{Cahit} the following theorem:
\begin{thm}\label{Cahit}
The following families of graphs are cordial:
\begin{itemize}
  \item [1.] trees;
  \item [2.] complete graphs $K_n$ if and only if $n \leqslant 3$;
  \item [3.] complete bipartite graphs $K_{n,m}$ for all $m$ and
$n$;
\item[4.] cycles $C_n$ if and only if $n \not \equiv 2 \mod
4$.
\end{itemize}

\end{thm}

However, Cairnie and  Edwards proved that in a general case the
problem of deciding whether or not a graph $G$ is cordial is
NP-complete \cite{CE}.\\

Hovey \cite{Hovey} introduced $k$-cordial labeling of graphs as a generalization
of harmonious and cordial labelings. He showed~\cite{Hovey} the following:
\begin{thm} \label{Hovey1}
All caterpillars are $k$-cordial for all $k$ and all trees are $k$-cordial for $k=3,4,5$.
\end{thm}

Moreover he advances~\cite{Hovey} the following conjecture:
\begin{cj} \label{Hovey2}
All trees are $k$-cordial for all $k$.
\end{cj}

\indent In this paper we consider the corresponding problem for
hypertrees. We show that $p$-uniform
hypertrees are $k$-cordial for some values of $k$.

\indent The paper is organized as follows. In the next section we
prove some preliminary lemmas. They will be needed in the proof of
the main theorem presented in the third section.
\section{Lemmas}
\begin{lem}\label{wspolnewierzch}
For any two edges $e_1 \neq e_2$ in a hypertree
there is at most one common vertex.
\end{lem}
\prf Suppose that $e_1\neq e_2$ have two vertices $v,u$ in
common. Thus the hypergraph has a cycle $v, e_1, u, e_2, v$.
\cbdo
\begin{lem}\label{rozmiar}
Let $T$ be a $p$-uniform hypertree. Then $|T|=(p-1)\|T\|+1$
\end{lem}
\prf The proof is by induction on the size of $T$. If $\|T\|=1$,
then the claim obviously holds. Assume the claim holds for every
$p$-uniform hypertree with size $k-1$. Let $T$ be a $p$-uniform
hypertree with size $k$. Let $e=x_1...x_p$ be the last edge in the
longest path in $T$ such that $d(x_1)=...=d(x_{p-1})=1$. Let
$T'=T-\{e\}$ be a hypertree with the vertex set $V'=V
\setminus\{x_1,...,x_{p-1}\}$ and the edge set $E'=E
\setminus\{e\}$. Observe, that $T'$ is a $p$-uniform hypertree. By
induction $$|T'|=(p-1)\|T'\|+1$$ hence
$$|T|=|T'|+p-1=(p-1)(\|T'\|+1)+1=(p-1)\|T\|+1.$$ \cbdo
\begin{lem}\label{sumy}
Let $k \geq 2$ be an integer. For every $a \in Z_k$ and every $1 \leq l<k$ there exist $x_1,...,x_l \in Z_k$ such that $x_i \neq x_j$ for $i \neq j$ and $x_1+...+x_l=a$.
\end{lem}
\prf
We distinguish two cases:\\
{\bf Case 1.} Let $k$ be an odd number. Then $\sum_{i=0}^{k-1}i=0\mod k$. Thus $Z_k$ contains $k$ mutually different elements such that their sum is equal to zero. Moreover, for every $x \in Z_K$ such that $x \neq 0$ the inverse element $k-x$ is other than $x$. Suppose first that $a=0$. For every $l\leq k$ let $A_l$ denotes a set with size $l$ of mutually different elements of $Z_k$ such that $\sum_{x \in A_l}x=0$. For odd $l<k$ we construct $A_l$ removing from $Z_k$ an adequate number of pairs $\{x,k-x\}$. For $l<k$ even we remove the zero and an adequate number of pairs $\{x,k-x\}$. Observe that for any $1<l \leq k$ and for any $a \in Z_k$ such that $a \neq 0$ we can construct $A_l$ this way that the pair $\{a, k-a\} \subset A_l$. \\
\indent Let $a \in Z_k$ such that $a \neq 0$. Denote by $B_l$ a set with size $l<k$ of mutually different elements of $Z_k$ such that $\sum_{x \in B_l}x=a$. We obtain $B_l$ removing from a proper set $A_{l+1}$ the element $k-a$.\\
{\bf Case 2.} Let $k$ be an even number. Then $\sum_{i=0}^{k-1}i=\frac{k}{2} \mod k$. Thus $Z_k$ contains $k$ mutually different elements such that their sum is equal to $\frac{k}{2}$. Observe that for every $x \in Z_k$ such that $x \neq 0$ and $x \neq \frac{k}{2}$ the element $\frac{k}{2}-x$ is other than $\frac{k}{2}$. Assume first that $a=\frac{k}{2}$. For every $l\leq k$ let $C_l$ denotes a set with size $l$ of mutually different elements of $Z_k$ such that $\sum_{x \in C_l}x=\frac{k}{2}$. For even $l<k$ we construct $C_l$ removing from $Z_k$ an adequate number of pairs $\{x,k-x\}$. For odd $l<k$ we remove the zero and an adequate number of pairs $\{x,k-x\}$ (thus $C_2=\{0, \frac{k}{2}\}$ and $C_1=\{\frac{k}{2}\}$). Observe that for any $2<l \leq k$ and for any $a \in Z_k$ such that $a \neq 0$ and $a \neq \frac{k}{2}$ we can construct $A_l$ this way that the pair $\{a, k-a\} \subset C_l$ and $\frac{k}{2} \in C_l$. \\
\indent Let $a \in Z_k$ such that $a \neq \frac{k}{2}$. Denote by $B_l$ a set with size $l<k$ of mutually different elements of $Z_k$ such that $\sum_{x \in B_l}x=a$. For $k>l>2$ we obtain $B_l$ removing from a proper set $C_{l+1}$ the element $\frac{k}{2}-a$ (in particular, the element $\frac{k}{2}$ if $a=0$). The existence of $B_l$ for $l=2$ or $l=1$ is obvious.  \cbdo
\section{The main result}
\begin{thm}\label{cordial}
Let $T$ be a $p$-uniform hypertree. $T$ is $k$-cordial if one of the following conditions hold:\\
$\bullet$ $p$ is odd and $k$ is even\\
$ \bullet$ $p=1 \mod k$\\
$\bullet$ $p=0 \mod k$.
\end{thm}
\prf For $k=1$ the above theorem is obvious so let us suppose that $k \geq 2$. The proof is by induction on the size of a hypertree. The above theorem obviously holds for any $p$-uniform hypertree with size one. Let $T$ be a $p$-uniform hypertree with size $\|T\|$ and assume that the theorem holds for every $p$-uniform hypertree with size less than $\|T\|$. Let $e=x_1...x_p$ be the last edge in the longest path in $T$ such that $d(x_1)=...=d(x_{p-1})=1$. Let $T'=T-\{e\}$ be a $p$-uniform hypertree with the vertex set $V'=V \setminus\{x_1,...,x_{p-1}\}$ and the edge set $E'=E \setminus\{e\}$. By induction there exists a $k$-friendly coloring $c'$ for $T'$ which induces $k$-cordial coloring ${c'}^*$. Below we show that we can extent the coloring $c'$ to a $k$-friendly coloring $c$ of $T$ this way that $c$ induces $k$-cordial coloring for $T$. So, we label vertices $x_1,...,x_{p-1}$ this way that we obtain a proper label $i_e=\sum_{i=1}^{p}c(x_i) \mod k$ of $e$. \\
\indent Let $|T'|=a \mod k$. Therefore there exists a subset $I \subset Z_k$ with size $a$ (in particular, $I$ is an empty set if $a=0$) such that $|v(c_i)|=|v(c_j)|$ for every $i,j \in I$ or $i,j \notin I$ and $|v(c_i)|=v(c_j)| +1$ for every $i \in I$, $j \notin I$. Then we label vertices $x_1,...,x_{k-a}$ using every element from $Z_k \setminus I$ exactly ones. Then, for every $i,j \in Z_k$ the number of vertices just labeled by $i$ is equal to the number of vertices labeled by $j$ and there are $p-1-(k-a)$ vertices not labeled in $T$. \\
\indent Let $r$ be an integer such that $p-1-(k-a)=r \mod k$. Observe that if $p=1 \mod k$ or $p$ is odd and $k$ is even, then $r>0$. We distinguish two cases:\\
{\it Case 1.} $r>0$. Then we color vertices not labeled before in the following way. Let $b$ be an integer such that $p-1-(k-a)=bk+r$. We part the set of vertices not labeled before into $k+1$ disjoint sets $U_0,...,U_{k-1}, R$ where each $U_i$ has size $b$ (in particular, $U_i$ is an empty set for $b=0$) and $R$ has size $r$. We color vertices from $U_i$ by $i$. So, for every $i,j \in Z_k$ the number of vertices just labeled by $i$ is equal to the number of vertices labeled by $j$ and there are $r$ vertices not labeled in $T$. Then, by Lemma~\ref{sumy} we can color $r$ remaining vertices from $R$ by $r$ different colors in such a way that we obtain an adequate color $i_e$ of $e$.\\
{\it Case 2.} $r=0$. Then $|T|=(p-1)\|T\|+1=0 \mod k$. For $p=0 \mod
k$ we obtain that $\|T\|=1 \mod k$ and hence $\|T'\|=0 \mod k$. So,
let us label vertices in $T$ not labeled before this way that we
obtain a $k$-friendly coloring $c$ of $T$. Then $c$ induces
$k$-cordial coloring for $T$. \cbdo

The following corollary follows from the Theorem~\ref{cordial}.
\begin{col}
All $p$-uniform hypertrees are cordial.
\end{col}
There are several open problem that one can consider. The main is the following generalization of Hovey conjecture:
\begin{cj}
All $p$-uniform hypertrees are $k$-cordial for all $k$.
\end{cj}
Besides the results presented in this article, we also obtained partial results concerning cases $k=3,4,5$, but they do not seem to give a proof (or a counterexample) above conjecture.

\end{document}